\documentclass[10pt,a4paper]{article}
\usepackage{amsfonts,amssymb,amsmath,amsthm,graphicx}
\usepackage[english]{babel}
\newtheorem{thm}{Theorem}
\newcounter{punkt}
\newcommand{\para}{\par\vspace{3ex plus0.2ex minus0.2ex}%
\noindent\refstepcounter{punkt}\textbf{\arabic{punkt}.} }
\newcounter{abbildung}
\newenvironment{meinproof}%[1]%
    {\begin{trivlist} \item \emph{Proof.}} %#1\/:}}%
    {{}\hfill $\square$ \end{trivlist}} %%Hans {}\hfill erg\"{a}nzt und \/ gel\"{o}scht

\newcommand{\EE}{{\mathbb E}}
\newcommand{\RR}{{\mathbb R}}
\newcommand{\cH}{{\mathcal H}}
\newcommand{\cQ}{{\mathcal Q}}
\newcommand{\cS}{{\mathcal S}}
\newcommand{\cT}{{\mathcal T}}
\newcommand{\T}{{\mathrm{T}}}
\newcommand{\tr}{{\mathrm{tr\,}}}
\newcommand{\rank}{{\mathrm{rank\,}}}
\let\phi=\varphi

\let\theta=\vartheta

\newcommand{\DelimArray}[4]{\left#1\begin{array}{*{#3}{c}}#4\end{array}\right#2}

\newcommand{\Mat}{\DelimArray()}

\newcommand{\pfeil}[1]{\overline{#1}\hspace{-0.6ex}\vec{\rule{0ex}{1.5ex}}\,}
\newcommand{\tetra}{{T}}
\newcommand{\derselbe}{{--------}}
\date{}
\begin{document}
\title{Altitudes of a Tetrahedron\\ and Traceless Quadratic Forms}
\author{Hans Havlicek \and Gunter Wei{\ss}}

\maketitle

 \para\textbf{INTRODUCTION.}
 \label{para:einleitung}
It is well known that the three altitudes of a triangle are
concurrent at the so-called orthocenter of the triangle. So one might
expect that the altitudes of a tetrahedron also meet at a point.
However, it was already pointed out in 1827 by the Swiss geometer
Jakob Steiner (1796--1863) that the altitudes of a general
tetrahedron are mutually skew, for they are generators of an
equilateral hyperboloid. This is a hyperboloid with the following
rather peculiar property: each nontangential plane that is
perpendicular to a generator meets the hyperboloid along an
equilateral hyperbola, i.e., a hyperbola with orthogonal asymptotes.
 \par
There are many papers, especially from the nineteenth century, that
deal with the altitudes of tetrahedra and related topics. An
excellent survey article with many historical footnotes is the paper
of N.A.\ Court \cite{court-48}. Another major source is the article
of M.\ Zacharias in the \emph{Encyclopedia of Mathematical Sciences}
\cite{zach-13} (completed in 1913). Also, the book of H.\ Schr\"oter
\cite{schroe-80} contains in section 28 some interesting remarks on
older papers that are not cited elsewhere.
 \par
Some of the many special points and lines associated with triangles
extend to tetrahedra, others do not. More precisely, there are often
a number of possibilities for generalizing a concept (like the
orthocenter) from the plane to $3$-space (or even $n$-space), since
separate notions in higher-dimensional geometry may coincide in the
plane. We illustrate this in section \ref{para:monge}, where we
discuss the Monge point of a tetrahedron. If the definition of the
Monge point is applied directly to the planar case, then the
orthocenter of a triangle is obtained. Table 1 starts with some
well-known ``noteworthy points'' of a triangle and lists their
existence in higher dimensions, together with short remarks or
references.
\par
%%\begin{table}[h]
\begin{center}{\small
\begin{tabular}{|l|c|c|c|c|}
 \hline
 Dimension     & $n=2$          & $n=3$        & $n\geq 4$        \\
 \hline
 Object        & Triangle       & Tetrahedron                   & $n$-Simplex \\
 \hline\hline
 Centroid $G$  & Yes            & Yes                           & Yes (See \cite{mehmke-84}, \cite{rose+jagl-69}.) \\
  \hline
 Circumcenter $C$ & Yes         & Yes                           & Yes (See \cite{mehmke-84}, \cite{rose+jagl-69}.) \\
  \hline
 Incenter $I$  & Yes            & Yes                           & Yes (See \cite{mehmke-84}, \cite{rose+jagl-69}.) \\
  \hline
 Orthocenter $H$& Yes           & Sometimes                     & Sometimes   \\
               &                & (sec.\ \ref{para:hoehen})    & (See \cite{mandan-62a}, \cite{mehmke-84}.)\\
  \hline
 Monge point $M$ & Yes: $M=H$   & Yes (sec.\ \ref{para:monge}) & Yes (See \cite{mandan-62a}, \cite{mehmke-84}.) \\
    &  (sec.\ \ref{para:monge})&                               & \\
  \hline
 Euler line $e$&Yes             & Yes (sec.\ \ref{para:euler}) & Yes (See \cite{mandan-62a}, \cite{mehmke-84}.)  \\
 ($G,C,M\in e$)&($G,C,H\in e$)  &                               & \\
 \hline
\end{tabular}
}%%end small
\\[1ex]\centerline{Table 1.}
\end{center}
%%\end{table}
 \par
Clearly, this list is far from being complete. For example,
Feuerbach's famous \emph{$9$-point circle} for the triangle is
generalized in \cite{fritsch-76} to the $n$-dimensional case as a
\emph{Feuerbach hypersphere}, and in \cite{court-48} two
\emph{$12$-point spheres} for an orthocentric tetrahedron are
described. However, in higher dimensions many questions are still
open and have not been treated systematically. Further references
are given in section \ref{para:ende}.
 \par
There are two major reasons for revisiting the subject ``altitudes of
a tetrahedron.'' On the one hand, we want to visualize the results,
since there are hardly any pictures in the cited papers. The figures
in this article have been prepared with the computer algebra system
Maple. Clearly, a figure must not replace a formal proof, but we are
convinced that figures can assist the reader in better understanding
spatial relationships. Also, in certain instances the idea behind a
proof is rather immediate from an adequate illustration.
 \par
On the other hand, we aim at a modern coordinate-free presentation in
terms of analytic geometry based on a Euclidean vector space, whereas
a lot of papers on the subject are written in terms of synthetic
projective geometry. Thus the prerequisite for reading this article
is knowledge of basic linear algebra.
 \par
We shall see that the quadric surface carrying the altitudes of a
general tetrahedron appears as a level set of a traceless (trace
$=0$) quadratic form $Q^*$. In fact, there is a natural link between
tetrahedra and certain traceless quadratic forms that will lead us to
an explicit expression for $Q^*$.

 \para\textbf{ALTITUDES AND ORTHOCENTRIC PERPENDICULARS.}
 \label{para:hoehen}
Let $\{A_0,A_1,A_2,A_3\}$ be the set of vertices of a nondegenerate
tetrahedron $\tetra$ in Euclidean $3$-space, i.e., assume that these
four points are not coplanar. Throughout this paper it is tacitly
assumed that $i$, $j$, $k$, and $l$ are indices subject to the
requirement that $\{i,j,k,l\}=\{0,1,2,3\}$. So, if $l$ is the index
of the ``top vertex,'' then $i$, $j$, and $k$ are the indices of the
``base.'' Since we are going to describe points by vectors of an
underlying Euclidean vector space $\EE$, we have to fix an origin
$O$. For the moment $O$ is chosen arbitrarily, although in section
\ref{para:quadrik} we make a specific choice of origin in order to
simplify our calculations.
 \par
Each vertex $A_i$ is given by its position vector
$a_i:=\pfeil{OA_i}$. Also, we introduce the vectors
\begin{equation}\label{eq:def-b_ij}
  b_{ij}:=a_i-a_j.
\end{equation}
These vectors have the following basic properties:
\begin{eqnarray}
 \label{eq:b-antisymm}
 &  b_{ij} + b_{ji}=0,& \\
 \label{eq:b-komplanar}
 & b_{ij} +b_{jk}+ b_{ki}=0,& \\
 \nonumber%%\label{eq:b-kantenzug}
 & b_{ij} +b_{jk}+ b_{kl}+b_{li}=0,& \\
 \label{eq:b-basis}
 & b_{ij},\ b_{ik},\ b_{il}\mbox{ are linearly independent}.&
\end{eqnarray}
Observe that (\ref{eq:b-basis}) follows from our assumption that
$\tetra$ is nondegenerate. The definition (\ref{eq:def-b_ij}) and
properties of the dot product in $\EE$ give the identity
\begin{equation}\label{eq:pluecker}
  b_{ij}\cdot b_{kl} + b_{ik}\cdot b_{lj}+b_{il}\cdot b_{jk}=0.
\end{equation}
Its geometric meaning will be explained in due course.
 \par
We consider an edge $A_iA_j$ of $\tetra$ and one of the remaining two
vertices, say $A_k$. The plane through $A_k$ perpendicular to the
specified edge is the set of all points $X$ described by the equation
\begin{equation}\label{eq:eckenebene}
  b_{ij}\cdot (a_k-x) = 0,
\end{equation}
where $x:=\pfeil{OX}$ denotes an ``unknown vector'' in $\EE$. We
observe that there are at most twelve such planes. They will be the
key to many of our considerations.
\par
We write $h_l$ for the altitude of $T$ passing through $A_l$. It lies
in every plane through $A_l$ that is perpendicular to the opposite
face. Hence $h_l$ lies in the three planes with equations
\begin{eqnarray}
  b_{ij}\cdot (a_l -x) &=&  0,\nonumber\\
  b_{jk}\cdot (a_l -x) &=&  0,\label{eq:h-lgs}\\
  b_{ki}\cdot (a_l -x) &=&  0.\nonumber
\end{eqnarray}
It follows from (\ref{eq:b-basis}) and (\ref{eq:b-antisymm}) that any
two of these equations are linearly independent, whereas
(\ref{eq:b-komplanar}) implies that all three equations are linearly
dependent. Thus (\ref{eq:h-lgs}) describes three different planes
that meet at the altitude $h_l$, any two of which are already
sufficient to determine $h_l$. Figure \ref{abb1} shows those planes
through the altitude $h_3$ that are perpendicular to $A_0A_1$ and
$A_0A_2$, respectively.
 \par
We now consider those uniquely determined planes through $A_k$,
$A_i$, and $A_j$ that are parallel to the first, the second, and the
third plane of (\ref{eq:h-lgs}), respectively. This yields the linear
system
\begin{eqnarray}
  b_{ij}\cdot (a_k -x) &=&  0,\nonumber\\
  b_{jk}\cdot (a_i -x) &=&  0,\label{eq:n-lgs}\\
  b_{ki}\cdot (a_j -x) &=&  0.\nonumber
\end{eqnarray}
Again, this system describes three planes that intersect in a line,
since
\begin{equation}\label{eq:n-bueschel}
  b_{ij}\cdot a_k + b_{jk}\cdot a_i + b_{ki}\cdot a_j =0
\end{equation}
follows from (\ref{eq:def-b_ij}), and $(b_{ij}+b_{jk}+b_{ki})\cdot
x=0\cdot x = 0$ is immediate from (\ref{eq:b-komplanar}). Let us
write $n_l$ for the line given by (\ref{eq:n-lgs}). (Two of these
planes through $n_3$ are drawn in Figure \ref{abb1}.) From its
definition, the line $n_l$ contains the orthocenter of the triangle
with vertices $A_i$, $A_j$, and $A_k$. Also, it is perpendicular to
the plane of that triangle. We call $n_l$ an \emph{orthocentric
perpendicular} of $\tetra$.
%%%%%%%%%%%%%%%%%%%%%%%%%%%%%%%%%%%%%%%%%%%%%%%%%%%%%%%%%%%%%%%%%%
{\unitlength1.3cm
%%   %\begin{figure}[ht]
      \begin{center}
      %% Figur Tetraeder mit Hoehe, Lot und Ebenen
      %% 1200dpi =  472,4409448819pixel/cm
      %% 3900x3900 = 4.13 x 4.13 cm bei 2400dpi
      \begin{minipage}[t]{6.669cm}%%{5.13cm}
         \begin{picture}(5.13,4.13)
         \put(0.5 ,0.0){\includegraphics[width=5.369cm]{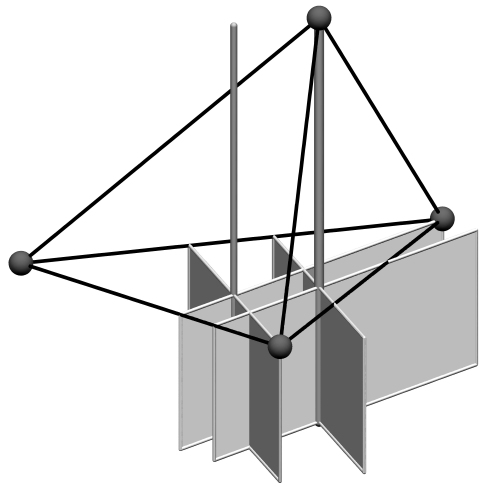}}
         \put(0.25,2.1){$A_0$}
         \put(2.87, 0.8){$A_1$}
         \put(4.4 ,2.4){$A_2$}
         \put(3.4 ,3.9){$A_3$}
         \put(3.3 ,2.6){$h_3$}
         \put(2.1 ,2.4){$n_3$}
         \end{picture}
         {\refstepcounter{abbildung}\label{abb1}
          \centerline{Figure \ref{abb1}.}}
      \end{minipage}
      \hspace{0cm}
      %%Figur Lot trifft Hoehen
      %% 1200dpi =  472,4409448819pixel/cm
      %% 3400x3400 = 3.60 x 3.60 cm bei 2400dpi
      %% vertikal angepasst!
      \begin{minipage}[t]{5.98cm}%%{4.6cm}
      %%vertikal:+0.265
         \begin{picture}(4.6,3.6)(0,-0.5)
         \put(0.5 ,0.0){\includegraphics[width=4.68cm]{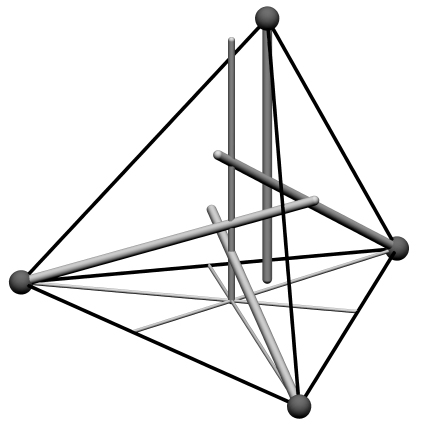}}
         \put(0.2 ,1.4){$A_0$}
         \put(3.2,0.1){$A_1$}
         \put(4.  ,1.6){$A_2$}
         \put(2.95,3.35){$A_3$}
         \put(2.12,2.5){$n_3$}
         \end{picture}
         {\refstepcounter{abbildung}\label{abb2}
          \centerline{Figure \ref{abb2}.}}
         \end{minipage}
      \end{center}
   %\end{figure}
}%
%%%%%%%%%%%%%%%%%%%%%%%%%%%%%%%%%%%%%%%%%%%%%%%%%%%%%%%%%%%%%%%%%%
The altitude $h_l$ and the perpendicular $n_l$ are parallel, but in
general they are not the same. However, we have the following result.
\begin{thm}\label{thm:1}
Each orthocentric perpendicular $n_l$ meets every altitude $h_i$ with
$i\neq l$.
\end{thm}
\begin{meinproof}
We read off from the second equation in (\ref{eq:h-lgs}) (after
interchanging $i$ with $l$) and the second equation in
(\ref{eq:n-lgs}) that $h_i$ and $n_l$ are both contained in the plane
with equation $b_{jk}\cdot (a_i-x)=0$. Furthermore, because $h_i$ and
$n_l$ are orthogonal to different faces of the tetrahedron $\tetra$,
they cannot be parallel. As a result, the two lines have a point in
common.
\end{meinproof}
\par
In Figure \ref{abb2} the line $n_3$ meets three out of the four
altitudes of $\tetra$. Also, under orthogonal projection onto the
plane $A_0A_1A_2$ the line $n_3$ is mapped to the orthocenter of the
triangle with vertices $A_0$, $A_1$, $A_2$, and the altitudes $h_0$,
$h_1$, and $h_2$ project to the altitudes of that triangle. We can
see that the foot of $h_3$ is on no altitude of the opposite
triangle. Hence $h_3$ does not meet any other altitude of the
tetrahedron. So we can seek a criterion for deciding whether or not
two altitudes meet at a point. It will be convenient to say that two
lines (with or without a common point) are \emph{orthogonal} if their
direction vectors are orthogonal.

\begin{thm}\label{thm:A}
An altitude $h_i$ meets an altitude $h_j$ precisely when
\begin{equation}\label{eq:schnittbed}
  b_{kl}\cdot b_{ij}=0
\end{equation}
or, in other words, when the opposite edges $A_kA_l$ and $A_iA_j$ are
orthogonal.
\end{thm}
\begin{meinproof}
From (\ref{eq:h-lgs}), the altitude $h_i$ is described by the linear
system
\begin{eqnarray*}
 b_{jk}\cdot (a_i-x)&=&0,\\
 b_{kl}\cdot (a_i-x)&=&0,
\end{eqnarray*}
whereas $h_j$ is determined by
\begin{eqnarray*}
 b_{kl}\cdot (a_j-x)&=&0,\\
 b_{li}\cdot (a_j-x)&=&0.
\end{eqnarray*}
Now (\ref{eq:schnittbed}) implies that $b_{kl}\cdot a_{i}=b_{kl}\cdot
a_{j}$. So two of these four equations are the same. Hence $h_i\cap
h_j\neq\emptyset$, for by (\ref{eq:b-basis}) the three planes given
by those equations are not parallel to a line.
\par
Conversely, let $P$ be a common point of $h_i$ and $h_j$. Put
$p=\pfeil{OP}$. We infer from
\begin{equation*}
  0=0-0=b_{kl}\cdot (a_i-p)-b_{kl}\cdot (a_j-p)=b_{kl}\cdot b_{ij}
\end{equation*}
that (\ref{eq:schnittbed}) is satisfied.
\end{meinproof}
\par
By symmetry, condition (\ref{eq:schnittbed}) is also necessary and
sufficient for the altitudes $h_k$ and $h_l$ to have a point in
common. Hence, if two altitudes intersect, then so do the other two.
Moreover, the reader will easily verify that (\ref{eq:schnittbed}),
$n_i\cap n_j\neq\emptyset$, and $n_k\cap n_l\neq\emptyset$ are
mutually equivalent.
 \par
At a first glance the following result may be somewhat surprising:
\begin{thm}\label{thm:B}
If one altitude meets two other altitudes then all altitudes are
concurrent.
\end{thm}
\begin{meinproof}
Suppose that $h_i$ meets $h_j$ and $h_k$. Then (\ref{eq:schnittbed})
implies $b_{kl}\cdot b_{ij}=b_{jl}\cdot b_{ik}=0$ and
(\ref{eq:pluecker}) shows that $b_{il}\cdot b_{jk}=0$, i.e., $h_i$
meets also $h_l$. By Theorem \ref{thm:A}, any two altitudes
intersect. Clearly, the four altitudes are not coplanar. So they have
a common point.
\end{meinproof}
\par
A tetrahedron with exactly one pair of opposite orthogonal edges is
called \emph{semi-orthocentric}. Figure \ref{abb3} shows such a
tetrahedron ($A_0A_3 \perp A_1A_2$). Observe that the foot of $h_3$
is on exactly one altitude of the opposite triangle. In Figure
\ref{abb4}, however, all four altitudes are concurrent at a point.
Such a tetrahedron is said to be \emph{orthocentric}, the
\emph{orthocenter} being the point where the altitudes meet. In an
orthocentric tetrahedron every edge is orthogonal to its opposite
edge. Furthermore, each altitude $h_i$ coincides with the
orthocentric perpendicular $n_i$.
\par
%%%%%%%%%%%%%%%%%%%%%%%%%%%%%%%%%%%%%%%%%%%%%%%%%%%%%%%%%%%%%%%%%%%%%%%%%%%%%%%
{\unitlength1.3cm
%%   %\begin{figure}[ht]
      \begin{center}
      %%Figur Hoehen schneiden paarweise
      %% 1200dpi =  472,4409448819pixel/cm
      %% 3600x3000 = 3,81 x 3,18 cm bei 2400dpi
      \begin{minipage}[t]{6.253cm}%%{4.81cm}
         \begin{picture}(4.81,3.18)
         \put(0.5,0.0)
         {\includegraphics[width=4.953cm]{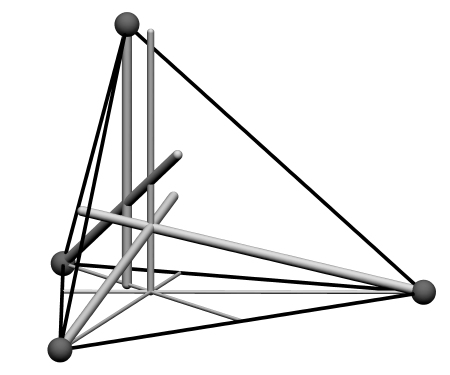}}
         \put(0.5 ,0.9) {$A_0$}
         \put(0.55 ,0.1){$A_1$}
         \put(4.25,0.6) {$A_2$}
         \put(1.05,2.85) {$A_3$}
         \put(1.85,2.15){$n_3$}
         \end{picture}
         {\refstepcounter{abbildung}\label{abb3}
          \centerline{Figure \ref{abb3}.}}
      \end{minipage}
      \hspace{0cm}
      %%Figur Hoehenschnittpunkt existiert
      %% 1200dpi =  472,4409448819pixel/cm
      %% 3600x3000 = 3,81 x 3,18 cm bei 2400dpi
      \begin{minipage}[t]{6.253cm}%%{4.81cm}
         \begin{picture}(4.81,3.18)
         \put(0.5,0.0)
         {\includegraphics[width=4.953cm]{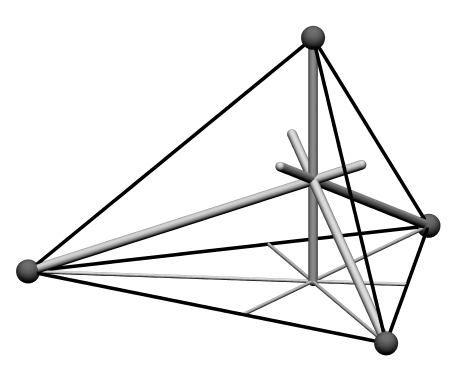}}
         \put(0.2 ,0.8){$A_0$}
         \put(3.9 ,0.15){$A_1$}
         \put(4.3 ,1.2){$A_2$}
         \put(3.3 ,2.8){$A_3$}
         \end{picture}
         {\refstepcounter{abbildung}\label{abb4}
          \centerline{Figure \ref{abb4}.}}
      \end{minipage}
      \end{center}
   %\end{figure}
}%
%%%%%%%%%%%%%%%%%%%%%%%%%%%%%%%%%%%%%%%%%%%%%%%%%%%%%%%%%%%%%%%%%%%%%%%%%%%%%%%
We refer to  \cite[p.~371]{bolt+jagl-69}, \cite{court-46}, and
\cite{court-64} for other proofs of Theorems \ref{thm:1},
\ref{thm:A}, and \ref{thm:B}.

\para\textbf{THE MONGE POINT OF A TETRAHEDRON.}\label{para:monge}
For an orthocentric tetrahedron the intersection of all planes given
by (\ref{eq:eckenebene}) is its orthocenter. However, we focus on an
arbitrary tetrahedron $\tetra$. Then among the planes
(\ref{eq:eckenebene}) there are two that are perpendicular to
$A_iA_j$ and pass through $A_k$ and $A_l$, respectively. These two
planes are either identical or disjoint. (In Figure \ref{abb1} two
such pairs of parallel planes can be seen.) In either case their
\emph{midplane} is the plane with equation
\begin{equation}\label{eq:mittenebene}
  b_{ij}\cdot(a_k+a_l-2x)=0.
\end{equation}
This midplane is orthogonal to the edge $A_iA_j$, and it passes
through the midpoint of the opposite edge $A_kA_l$ rather than the
midpoint of the edge $A_iA_j$, as the \emph{perpendicular bisector}
of $A_iA_j$ does. In general, the midplane (\ref{eq:mittenebene}) is
therfore \emph{not} the perpendicular bisector of the edge $A_iA_j$.
The equation of that plane is given in equation
(\ref{eq:symm_ebene}). The tetrahedron $\tetra$ has six midplanes.
 \par
Now for an arbitrary tetrahedron $\tetra$ there will be no
orthocenter, but $\tetra$ will have a point discovered by Gaspard
Monge (1746--1818) that is now called its \emph{Monge point}. The
construction of this point goes as follows:

\begin{thm}\label{thm:C}
All six midplanes of a tetrahedron are concurrent at a point.
\end{thm}
\begin{meinproof}
From (\ref{eq:b-basis}) we know that the vectors $b_{01}$, $b_{02}$,
and $b_{03}$ are linearly independent. Hence there exists a unique
common point, say $M$, of the (mutually nonparallel) planes whose
equations are
\begin{eqnarray}
  b_{01}\cdot(a_2+a_3-2x)&=&0,\nonumber\\
  b_{02}\cdot(a_1+a_3-2x)&=&0,\nonumber\\ %%\label{eq:M_lgs}\\
  b_{03}\cdot(a_1+a_2-2x)&=&0.\nonumber
\end{eqnarray}
We subtract the first equation from the second and obtain
\begin{eqnarray*}
  0&=& b_{02}\cdot a_1-b_{01}\cdot a_2 +
  (b_{02}-b_{01})\cdot(a_3-2x)\\
  &=& b_{02}\cdot a_1+b_{10}\cdot a_2 +
  b_{12}\cdot(a_3-2x)\\
  &=& b_{12}\cdot(a_0+a_3-2x),
\end{eqnarray*}
where we have used equations (\ref{eq:b-antisymm}),
(\ref{eq:b-komplanar}), and (\ref{eq:n-bueschel}). So $M$ lies in the
midplane that is perpendicular to the edge $A_1A_2$. A similar
calculation shows that $M$ belongs to the remaining midplanes as
well.
\end{meinproof}
\par
Figure \ref{abb5} displays a tetrahedron and its Monge point $M$. It
lies on that line of the plane spanned by the parallel lines $h_i$
and $n_i$ which is equidistant from both; one such line (dotted) is
illustrated for $i=3$. If we deform $\tetra$ by ``sliding'' the
vertex $A_3$ down along the line $h_3$, holding all other vertices
fixed, then the center of gravity of $\tetra$ remains an inner point
of $\tetra$, while the circumcenter of $\tetra$ approaches ``negative
infinity,'' since the radius of the circumscribed sphere of $\tetra$
tends to infinity. We shall see in section \ref{para:euler} that the
center of gravity of $\tetra$ is the midpoint of $M$ and the
circumcenter of $\tetra$. Hence the Monge point $M$ ``moves up''
along the dotted line. This gives Figure \ref{abb6}, in which the
Monge point is exterior to the tetrahedron. For the proper choice of
$A_3$ we could also produce a tetrahedron whose the Monge point is
incident with a face.
%%%%%%%%%%%%%%%%%%%%%%%%%%%%%%%%%%%%%%%%%%%%%%%%%%%%%%%%%%%%%%%%%%
{\unitlength1.3cm
%%   %\begin{figure}[ht]
      \begin{center}
      %% Figur steiles Tetraeder mit Mongepunkt Hoehe, Lot
      %% 1200dpi =  472,4409448819pixel/cm
      %% 3900x3100 = 4.13 x 3.28 cm bei 2400dpi
      \begin{minipage}[t]{6.669cm}%%{5.13cm}
         \begin{picture}(5.13,3.28)
         \put(0.5 ,0.0){\includegraphics[width=5.369cm]{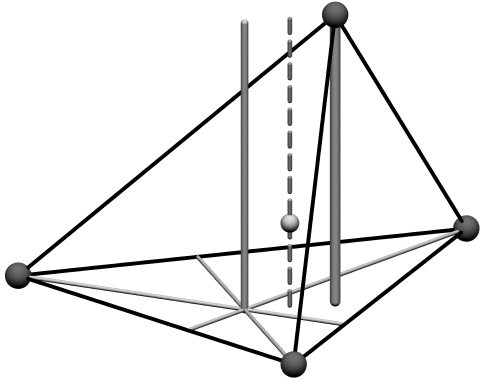}}
         \put(0.3 ,1.2){$A_0$}
         \put(3.07,0.0){$A_1$}
         \put(4.6 ,1.5){$A_2$}
         \put(3.5 ,3.1){$A_3$}
         \put(3.5 ,1.8){$h_3$}
         \put(2.2 ,1.6){$n_3$}
         \put(2.62,1.5){$M$}
         \end{picture}
         {\refstepcounter{abbildung}\label{abb5}
          \centerline{Figure \ref{abb5}.}}
      \end{minipage}
      \hspace{-0.7cm}
      %% Figur flaches Tetraeder mit Mongepunkt Hoehe, Lot
      %% 1200dpi =  472,4409448819pixel/cm
      %% 3900x3100 = 4.13 x 3.28 cm bei 2400dpi
      \begin{minipage}[t]{6.669cm}%%{5.13cm}
         \begin{picture}(5.13,3.28)
         \put(0.5 ,0.0){\includegraphics[width=5.369cm]{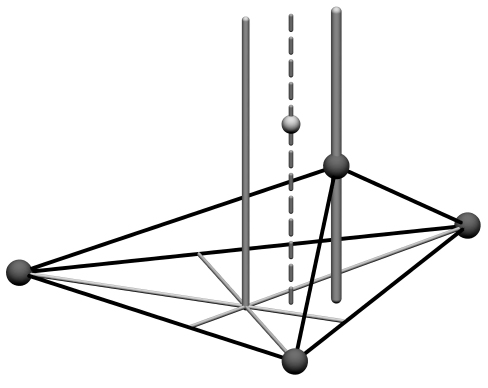}}
         \put(0.3 ,1.2){$A_0$}
         \put(3.07,0.0){$A_1$}
         \put(4.6 ,1.5){$A_2$}
         \put(3.5 ,1.9){$A_3$}
         \put(3.5 ,2.5){$h_3$}
         \put(2.2 ,1.8){$n_3$}
         \put(2.62,2.3){$M$}
         \end{picture}
         {\refstepcounter{abbildung}\label{abb6}
          \centerline{Figure \ref{abb6}.}}
         \end{minipage}
      \end{center}
   %\end{figure}
}%
%%%%%%%%%%%%%%%%%%%%%%%%%%%%%%%%%%%%%%%%%%%%%%%%%%%%%%%%%%%%%%%%%%
\par
Let $m:=\pfeil{OM}$ be the position vector of the Monge point $M$
described in Theorem \ref{thm:C}. Then (\ref{eq:mittenebene}) implies
that
   $b_{ij}\cdot(a_k+a_l-2m)= b_{kl}\cdot(a_i+a_j-2m)$ and
multiplying this out leads to the relation
\begin{equation*}
  2a_i\cdot a_l - 2a_i\cdot m -2a_l\cdot m
  = 2a_j\cdot a_k - 2a_k\cdot m -2a_j\cdot m,
\end{equation*}
which is equivalent to
\begin{equation}\label{eq:a-m}
  (a_i-m)\cdot (a_l-m) = (a_j-m)\cdot (a_k-m).
\end{equation}
This pretty equation illustrates in another way the special role of
the point $M$.
\par
Clearly, if all altitudes of $\tetra$ are concurrent, then the Monge
point $M$ is their point of intersection. This is one rationale for
the assertion that the Monge point serves as a substitute for the
``missing orthocenter'' of a general tetrahedron.
\par
We sketch the definition of the Monge point in a more general context
in order to illustrate that the orthocenter of a triangle can be
considered as its Monge point. An $n$-simplex $S$ in $n$-dimensional
Euclidean space has $n+1$ vertices. For each edge there is a unique
hyperplane that is perpendicular to that edge and that contains the
center of gravity of the $n-1$ ``opposite vertices,'' i.e., the
vertices not on the given edge. There are $n+1\choose 2$ such
hyperplanes, and they have a point in common---the Monge point of the
$n$-simplex (see \cite{mandan-62a} or \cite{mehmke-84}). For $n=3$
(tetrahedron) this is in accordance with Theorem \ref{thm:C}, since
each edge has two opposite vertices whose center of gravity is just
their midpoint. For $n=2$ (triangle) each edge has a single opposite
vertex that is its own center of gravity. In this situation, the
three altitudes of the triangle replace the hyperplanes from the
general case, and their common point (the orthocenter) coincides with
the Monge point.

 \para\textbf{THE EULER LINE IN SPACE.}\label{para:euler}
We consider the circumcenter $C$, the center of gravity $G$, and the
Monge point $M$ of the tetrahedron $\tetra $ with vertex set
$\{A_0,A_1,A_2,A_3\}$. An equation
\begin{equation}\label{eq:symm_ebene}
  b_{ij}\cdot(a_i+a_j-2x)=0
\end{equation}
describes the perpendicular bisector of the edge $A_iA_j$, so that
$\pfeil{OC}$ is the only solution of the linear system of all six
equations (\ref{eq:symm_ebene}).

If we fix indices $i$ and $j$, then the midplane of the planes
represented by (\ref{eq:mittenebene}) and (\ref{eq:symm_ebene}) has
an equation of the form
\begin{equation}\label{eq:schwerpunktebene}
  b_{ij}\cdot(a_i+a_j+a_k+a_l-4x)=0.
\end{equation}
This plane contains the midpoint of the segment $CM$. From
(\ref{eq:b-basis}) we know that there are three linearly independent
vectors among the $b_{ij}$, whence the only solution of the linear
system comprising all equations (\ref{eq:schwerpunktebene}) is
\begin{equation*}%%\label{eq:schwerpunkt}
  \pfeil{OG}=\frac{1}{4}(a_i+a_j+a_k+a_l).
\end{equation*}
This means that the center of gravity is the midpoint of the segment
$CM$. In other words, whenever two of the points $G$, $C$, and $M$
are different, their join can be considered as an analog of the Euler
line in $3$-space. (We remind the reader that the \emph{Euler line}
of a triangle contains the center of gravity, the circumcenter, and
orthocenter of the given triangle.)
 \par
We add in passing that in the $n$-dimensional setting there is also
an Euler line: the center of gravity divides the segment formed by
the circumcenter and the Monge point internally in the ratio
$2:(n-1)$ (see \cite{mandan-62a} or \cite{mehmke-84}). This explains
why in the plane (orthocenter = Monge point) the ratio on the Euler
line is $2:1$, whereas in three dimensions it is $2:2$.

\para\textbf{TRACELESS QUADRATIC FORMS.}\label{para:spurfrei}
If $Q:\EE\to\RR$ is a nonzero quadratic form and $\rho$ in $\RR$ is a
constant, then $Q(x)=\rho$ is an equation of a (possibly degenerate)
quadric surface that is symmetric with respect to the origin $O$. We
refer to \cite{gruen+w-77} for basic properties of quadrics and
quadratic forms.
\par
By the polarization formula, a quadratic form $Q:\EE\to\RR$ gives
rise to a symmetric bilinear form $f:\EE\times \EE\to\RR$,
\begin{equation*}%%\label{eq:polarform}
  f(v,w)=\frac{1}{2}\big(Q(v+w)-Q(v)-Q(w)\big),
\end{equation*}
with the property $Q(v)=f(v,v)$ for all vectors $v$ in $\EE$. This
function $f$ is called the polar form of $Q$.

Let $\{e_1,e_2,e_3\}$ be an orthonormal basis of $\EE$. Then each
vector $x$ is uniquely determined by its Cartesian coordinates
$(\xi_1,\xi_2,\xi_3)$ in $\RR^3$, viz.,
$x=\xi_1e_1+\xi_2e_2+\xi_3e_3$. The symmetric $3$-by-$3$ matrix
$\Sigma=(\sigma_{rs}):=(f(e_r,e_s))$ allows us to express $Q$ in
terms of coordinates in the form
\begin{equation*}%%\label{eq:Q_koo}
  Q(x) = \sum_{r,s=1}^{3} \sigma_{rs}\xi_r\xi_s.
\end{equation*}
If we change to another orthonormal basis, then $\Sigma$ changes to a
congruent matrix $\Sigma'=\Omega^\T\Sigma\Omega$, where $\Omega$ is
an orthogonal matrix and $\Omega^\T$ denotes the transpose of
$\Omega$. As $\Omega$ is orthogonal, we have $\Omega^\T=\Omega^{-1}$.
The matrices $\Sigma$ and $\Sigma'$ are thus similar. It is well
known that similar matrices have the same trace. So it makes sense to
speak of the \emph{trace} $\tr Q$ \emph{of a quadratic form} $Q$, as
long as we restrict ourselves to orthonormal bases.
 \par
Consider, for example, arbitrary vectors $c$ and $d$\ of $\EE$ with
Cartesian coordinates $(\gamma_1,\gamma_2,\gamma_3)$ and
$(\delta_1,\delta_2,\delta_3)$, respectively, and the quadratic form
\begin{equation}\label{eq:c.d}
  x\mapsto (x\cdot c)(x\cdot d).
\end{equation}
The $(r,s)$-entry of its matrix equals
$\frac{1}{2}(\gamma_r\delta_s+\gamma_s\delta_r)$, whence the trace of
(\ref{eq:c.d}) is
\begin{equation}\label{eq:tr_c.d}
 \gamma_1\delta_1+\gamma_2\delta_2+\gamma_3\delta_3=c\cdot d.
\end{equation}
 \par
We are particularly interested in quadratic forms with trace zero or,
in other words, \emph{traceless} quadratic forms. If $Q$ is an
arbitrary quadratic form on $\EE$ then, by transformation on
principal axes, there exists an orthonormal basis $\{e_1,e_2,e_3\}$
of $\EE$ with respect to which
\begin{equation}\label{eq:diagonal}
     Q(x)=\sigma_{11}\xi_1^2+\sigma_{22}\xi_2^2+\sigma_{33}\xi_3^2.
\end{equation}
Suppose now that $\tr Q=0$. Hence
\begin{equation}\label{eq:spurfrei}
   \tr Q = \sigma_{11}+\sigma_{22}+\sigma_{33}=0
\end{equation}
and there are several cases to consider. These depend on the rank of
$Q$, meaning the rank of any associated matrix.
 \par
Case 1: $\rank Q\leq 1$. Then at most one coefficient $\sigma_{rr}$
is nonzero and (\ref{eq:spurfrei}) shows that
$\sigma_{11}=\sigma_{22}=\sigma_{33}=0$, i.e., $Q$ is the zero-form.
 \par
Case 2: $\rank Q=2$. Then the basis can be chosen in such a way that
$\sigma_{33}=0$. This means that $\sigma_{22}=-\sigma_{11}\neq 0$ and
$Q(x)=\sigma_{11}(\xi_1-\xi_2)(\xi_1+\xi_2)$. Consequently, $Q(x)=0$
describes a \emph{pair of orthogonal planes} with equations
$\xi_1=\xi_2$ and $\xi_1=-\xi_2$, respectively.
 \par
Case 3: $\rank Q=3$. Then $Q$ is indefinite by (\ref{eq:spurfrei})
and $Q(x)=0$ is the equation of a quadratic cone (with vertex at the
origin) that is called an \emph{equilateral cone}. Its generators
(i.e., the lines entirely contained in the cone) have the following
remarkable property (see \cite[p.~293, Ex.\ 19]{semp+k-98}):
\begin{thm}\label{thm:gls-kegel}
If $g$ is a generator of an equilateral cone, then there are
generators $g_1$ and $g_2$ of the cone such that $g$, $g_1$, and
$g_2$ are mutually orthogonal.
\end{thm}
\begin{meinproof}
As earlier, let the cone be given in the form $Q(x)=0$, where $Q$ is
a traceless quadratic form of rank $3$. Choose $e_3$ to be a unit
vector in the direction of the generator $g$. Also, let $\{e_1,e_2\}$
be an orthonormal basis of $g^\perp$ (the orthogonal complement of
$g$) giving the principal axes of $Q$ restricted to $g^\perp$. Then
the matrix of the associated polar form $f$ reads
\begin{equation*}
  \Mat3{\sigma_{11} & 0            & \sigma_{13}\\
        0            & \sigma_{22} & \sigma_{23}\\
        \sigma_{31} & \sigma_{32} & 0 }
\end{equation*}
and $\rank Q=3$ implies that $\sigma_{11}=-\sigma_{22}\neq 0$. We
infer that the plane $g^\perp$ meets the cone in the orthogonal lines
$g_1:=\RR(e_1+e_2)$ and $g_2:=\RR(e_1-e_2)$.
\end{meinproof}
\par
We conclude that an equilateral cone carries infinitely many
orthogonal tripods. Figure \ref{abb7} shows an equilateral cone, some
of its orthogonal tripods (one in white, the others in black), and
the principal axes.
%%%%%%%%%%%%%%%%%%%%%%%%%%%%%%%%%%%%%%%%%%%%%%%%%%%%%%%%%%%%%%%%%%%%%%%%%%%%%%%
{\unitlength1.3cm
%%   %\begin{figure}[ht]
      \begin{center}
      %%Figur Gleichseitiger Kegel
      %% 1200dpi =  472,4409448819pixel/cm
      %% 3000x3300 = 3,18 x 3,49cm bei 2400dpi
      \begin{minipage}[t]{4.134cm}%%{3.18cm}
         \begin{picture}(4.18,3.49)
         \put(0.0,0.0)%%kein hor. Abstand
         {\includegraphics[width=4.134cm]{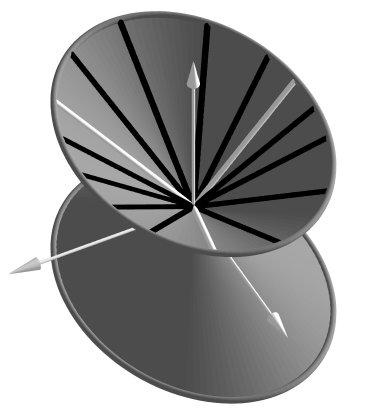}}
         \end{picture}
         {\refstepcounter{abbildung}\label{abb7}
          \centerline{Figure \ref{abb7}.}}
      \end{minipage}
      \hspace{1cm}
      %%Figur Dreibein
      %% 1200dpi =  472,4409448819pixel/cm
      %% 4800x2600 = 5,08 x 2,75 cm bei 2400dpi
      \begin{minipage}[t]{6.604cm}%%{5.08cm}
         \begin{picture}(5.08,3.75)
         \put(0.0,0.5)%%Vertikale Korrtur: Platz fuer $O$.
         {\includegraphics[width=6.604cm]{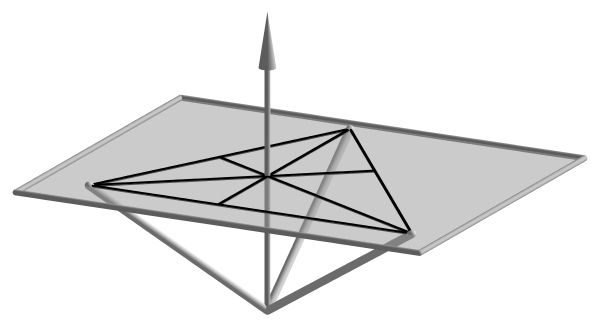}}
         \put(2.1 ,0.3){$O$}
         \put(0.75,1.9){$L_1$}
         \put(3.5 ,1.45){$L_2$}
         \put(2.9 ,2.4){$L_3$}
         \end{picture}
         {\refstepcounter{abbildung}\label{abb8}
          \centerline{Figure \ref{abb8}.}}
      \end{minipage}
      \end{center}
   %\end{figure}
}%
%%%%%%%%%%%%%%%%%%%%%%%%%%%%%%%%%%%%%%%%%%%%%%%%%%%%%%%%%%%%%%%%%%%%%%%%%%%%%%%
Theorem \ref{thm:gls-kegel} can be transferred to a result in the
plane as follows: Choose any orthogonal tripod on the cone and a
plane that meets the legs of the tripod in distinct points, say
$L_1$,$L_2$, and $L_3$ (see Figure \ref{abb8}). Hence
$\{O,L_1,L_2,L_3\}$ determines a so-called {trirectangular
tetrahedron}. The edges $OL_1$, $OL_2$, and $OL_3$ are at the same
time altitudes of this tetrahedron, and $O$ is its orthocenter.
Therefore the orthocentric perpendicular of the triangle $\Delta$
with vertices $L_1$, $L_2$, and $L_3$ runs through $O$ or, said
differently, the orthogonal projection of $O$ onto the plane
$L_1L_2L_3$ is just the orthocenter of this triangle. It is also
worth noting that $\Delta$ is always an acute triangle. Consider, for
example, the angle $\phi_1=\angle\, L_2L_1L_3$. From the law of
Pythagoras we get
\begin{equation*}
  \overline{L_1L_2}^2+\overline{L_1L_3}^2
  =\overline{OL_1}^2+\overline{OL_2}^2+\overline{OL_1}^2+\overline{OL_3}^2
  =\overline{L_2L_3}^2+2\overline{OL_1}^2.
\end{equation*}
Now the cosine law, applied to $\Delta$ and $\phi_1$, implies that
\begin{equation*}
  2\, \overline{OL_1}^2=2\, \overline{L_1L_2}\;\overline{L_1L_3}\,\cos\phi_1.
\end{equation*}
Since $\overline{OL_1}$, $\overline{L_1L_2}$, and $\overline{L_1L_3}$
are positive real numbers, so is $\cos\phi_1$, i.e., we have an acute
angle at $L_1$.
 \par
Let us turn back to the orthogonal tripods of an equilateral cone $C$
with an equation $Q(x)=0$ given in coordinates according to
(\ref{eq:diagonal}) and (\ref{eq:spurfrei}); we may assume
$\sigma_{11},\sigma_{22}>0$ and $\sigma_{33}<0$. We select a plane
$P: \xi_3=\rho$, where $\rho\neq 0$. Then $P$ meets $C$ along an
ellipse $E$. This gives Figure \ref{abb9}: the infinitely many
orthogonal tripods on the cone give rise to an infinite family of
acute triangles that are inscribed in $E$. Furthermore, the center of
$E$ is the common orthocenter of all these triangles. Clearly, if the
ellipse happens to be a circle ($\sigma_{11}=\sigma_{22})$, then all
triangles will be equilateral.
%%%%%%%%%%%%%%%%%%%%%%%%%%%%%%%%%%%%%%%%%%%%%%%%%%%%%%%%%%%%%%%%%%%%%%%%%%%%%%%
{\unitlength1cm
%%   %\begin{figure}[ht]
      \begin{center}
      %%Figur Schliessungskegelschnitt
      %% 1200dpi =  472,4409448819pixel/cm
      %% 4800x3500 = 5,08 x 3,70 cm bei 2400dpi
      \begin{minipage}[t]{6.08cm}
         \begin{picture}(6.08,3.70)
         \put(0.5,0.0)
         {\includegraphics[width=5.08cm]{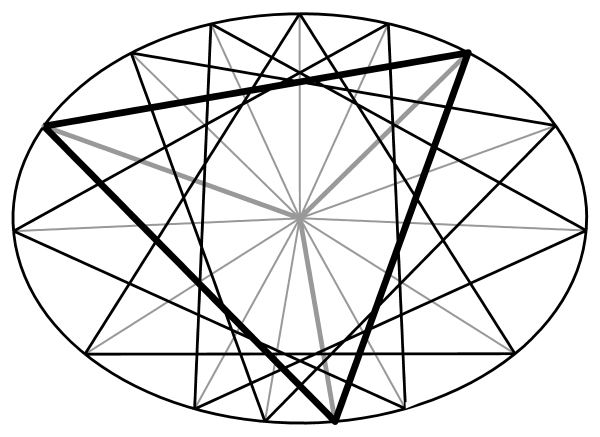}}
         \put(5.1,3){$E$}
         \end{picture}
         {\refstepcounter{abbildung}\label{abb9}
          \centerline{Figure \ref{abb9}.}}
         \end{minipage}
      \end{center}
   %\end{figure}
}%
%%%%%%%%%%%%%%%%%%%%%%%%%%%%%%%%%%%%%%%%%%%%%%%%%%%%%%%%%%%%%%%%%%%%%%%%%%%%%%%
Conversely, it is easily seen that every ellipse arises in this way
as a planar section of an equilateral cone. Hence the problem of
finding all triangles with the aforementioned properties always has
infinitely many solutions. Such a problem is called a \emph{porism}.
The investigation of problems of this kind goes back to the French
geometer J.V.\ Poncelet (1788--1867), and it is best understood in
terms of projective geometry: the triangles under consideration are
not only inscribed in the given ellipse but are also self-polar with
respect to an elliptic polarity. It follows that the sides of the
triangles are the tangents to another ellipse. The interested reader
should consult \cite[chaps.\ 4, 5]{semp+k-98}.

\para\textbf{THE QUADRIC OF THE ALTITUDES.}\label{para:quadrik}
In this section we choose the Monge point $M$ of a given tetrahedron
$\tetra$ as the origin. This will simplify our calculations.
Furthermore, we introduce the scalars
\begin{equation}\label{eq:def_lambda}
\lambda_{ij}:=a_i\cdot a_j  \; (i\neq j).
\end{equation}
From the symmetry of the dot product we have
$\lambda_{ij}=\lambda_{ji}$, but it is more important to notice that
from $m=\pfeil{MM}=0$ and (\ref{eq:a-m}) follows
\begin{equation}\label{eq:lambda_symm}
 \lambda_{il}  = \lambda_{jk};
\end{equation}
i.e., each $\lambda_{ij}$ is equal to one of the numbers
$\lambda_{01}$, $\lambda_{02}$, or $\lambda_{03}$.
 \par
First, let us rewrite the condition from (\ref{eq:schnittbed}) from
Theorem \ref{thm:A}. We see that the altitudes $h_i$ and $h_j$ meet
at a point precisely when
\begin{equation}\label{eq:schnittbed_lambda}
  b_{kl}\cdot b_{ij}=(a_k-a_l)\cdot (a_i-a_j)= 2
(\lambda_{ki}-\lambda_{li})=2 (\lambda_{lj}-\lambda_{kj})=0.
\end{equation}
Thus for a ``generic'' tetrahedron the scalars $\lambda_{01}$,
$\lambda_{02}$, and $\lambda_{03}$ are distinct, for a
semi-orthocentric tetrahedron exactly two of them are identical, and
for a tetrahedron to be orthocentric it is necessary and sufficient
that $\lambda_{01}=\lambda_{02}=\lambda_{03}$.
  \par
Next, we look for quadratic forms on $\EE$ that arise from $\tetra$
in a natural way. We start with the forms
\begin{equation*}%%\label{def:Q_ij}
Q_{ijkl}: x\mapsto (x\cdot b_{ij})(x\cdot b_{kl}).
\end{equation*}
Since $M$ has been chosen as the origin, every quadric
$Q_{ijkl}(x)=0$ is the union of two midplanes, namely, the planes
with equations $b_{ij}\cdot x=0$ and $b_{kl}\cdot x=0$ (cf.\
(\ref{eq:mittenebene})).
 \par
The collection of all quadratic forms $\EE\to\RR$ constitutes a real
vector space $\cQ$ isomorphic to the six-dimensional space of
symmetric $3$-by-$3$ matrices over $\RR$. However, we are only
interested in certain subspaces of $\cQ$:
\begin{thm}\label{thm:S_T}
 The subspace $\cS$ of $\cQ$ that is spanned by all quadratic forms $Q_{ijkl}$
has dimension two. The subspace $\cT$ of all traceless quadratic
forms in $\cS$ is either one- or two-dimensional.
\end{thm}
\begin{meinproof}
By (\ref{eq:b-antisymm}), $Q_{ijkl}=Q_{klij}=-Q_{jikl}$. So each
$Q_{ijkl}$ is equal to $\pm Q_{0123}$, $\pm Q_{0231}$, or $\pm
Q_{0312}$. From (\ref{eq:def-b_ij}) and a straightforward calculation
it follows that
\begin{equation}\label{eq:Q_ij}
  Q_{0123}(x)+Q_{0231}(x)+Q_{0312}(x)=0
\end{equation}
for all $x$ in $\EE$. We infer that $\cS$ is generated by $Q_{0123}$
and $Q_{0231}$. The midplanes perpendicular to $b_{01}$, $b_{23}$,
$b_{02}$, and $b_{13}$ are distinct. Hence we can find a vector $v$
with $Q_{0123}(v)\neq Q_{0231}(v)=0$. Since $Q_{0123}$ is not
proportional to $Q_{0231}$, $\dim\cS=2$.
 \par
The trace function is a nonzero linear form from $\cQ$ to $\RR$; its
restriction to $\cS$ is either nonzero, in which case $\dim\cT=1$, or
zero, in which event $\cT=\cS$ is two-dimensional.
\end{meinproof}
\par
From (\ref{eq:tr_c.d}) and (\ref{eq:schnittbed_lambda}) it follows
that
\begin{equation}\label{eq:tr_Q_ij}
 \tr Q_{ijkl} = b_{ij}\cdot b_{kl}
            = 2 (\lambda_{ki}-\lambda_{li})
            = 2 (\lambda_{lj}-\lambda_{kj}).
\end{equation}
We note in passing that (\ref{eq:Q_ij}) and (\ref{eq:tr_Q_ij})
illustrate the meaning of equation (\ref{eq:pluecker}): it simply
says that the zero form $Q_{0123}+Q_{0231}+Q_{0312}$ has trace $0$.
In addition, (\ref{eq:tr_Q_ij}) implies that $\tr Q_{ijkl}=0$ is
necessary and sufficient for the altitudes $h_i$ and $h_j$ to meet.
 \par
In particular, let us consider the quadratic form
\begin{equation}\label{eq:Q^*_def}
 Q^*:=
  \lambda_{01}Q_{0123}
 +\lambda_{02}Q_{0231}
 +\lambda_{03}Q_{0312}.
\end{equation}
Equations (\ref{eq:tr_Q_ij}) and the symmetry conditions
(\ref{eq:lambda_symm}) yield
\begin{equation}\label{eq:spurnull}
  \tr Q^* =
          \lambda_{01}(2(\lambda_{02}-\lambda_{03}))
        + \lambda_{02}(2(\lambda_{03}-\lambda_{01}))
        + \lambda_{03}(2(\lambda_{01}-\lambda_{02}))
        = 0.
\end{equation}
We are now in a position to state our main results: Theorem
\ref{thm:gleichung} states that the points on the altitudes satisfy
the quadratic equation (\ref{eq:gleichung}) if the Monge point $M$ is
chosen as the origin. However, it remains open here whether or not
$Q^*$ is identically zero. Thus equation (\ref{eq:gleichung}) may be
trivial ($0=0$). We shall see later that this is the case if and only
if the tetrahedron $\tetra$ is orthocentric. Otherwise,
(\ref{eq:gleichung}) is the equation of a (possibly degenerate)
quadric carrying the four altitudes. This will be shown in Theorem
\ref{thm:hauptsatz} and the subsequent remarks.

\begin{thm}\label{thm:gleichung}
Let $P$ be a point on an altitude of a tetrahedron with Monge point
$M$. Then $p:=\pfeil{MP}$ is a solution of the quadratic equation
\begin{equation}\label{eq:gleichung}
  Q^*(x) =
   (\lambda_{01}-\lambda_{02})(\lambda_{02}-\lambda_{03})(\lambda_{03}-\lambda_{01}),
\end{equation}
where the quadratic form $Q^*$ and the scalars $\lambda_{ij}$ are
given by (\ref{eq:Q^*_def}) and (\ref{eq:def_lambda}), respectively.
\end{thm}
\begin{meinproof}
If we apply a transposition on $(0,1,2,3)$, then both sides of
(\ref{eq:gleichung}) are multiplied by $-1$, where we have to take
into account the symmetry conditions (\ref{eq:lambda_symm}). So
equation (\ref{eq:gleichung}) remains unchanged, up to multiplication
by $\pm 1$, under any permutation of $(0,1,2,3)$. Therefore it is
enough to establish that all $p$ determined by points $P$ on the
altitude $h_0$ satisfy (\ref{eq:gleichung}).
\par
Formula (\ref{eq:Q_ij}) shows that
\begin{equation}\label{eq:Q^*_asymm}
 Q^*=-(\lambda_{03}-\lambda_{01})Q_{0123}+(\lambda_{02}-\lambda_{03})Q_{0231}.
\end{equation}
Furthermore, (\ref{eq:h-lgs}) implies that $p\cdot b_{ij}=a_0\cdot
b_{ij}=\lambda_{0i}-\lambda_{0j}$ whenever $i,j\neq 0$. Hence
\begin{eqnarray*}
\renewcommand\arraystretch{1.2}
 Q^*(p)& = & -(\lambda_{03}-\lambda_{01})(p\cdot b_{01})(p\cdot b_{23})
          +(\lambda_{02}-\lambda_{03})(p\cdot b_{02})(p\cdot b_{31})\\
       & = & -(\lambda_{03}-\lambda_{01})(p\cdot b_{01})(\lambda_{02}-\lambda_{03})
          +(\lambda_{02}-\lambda_{03})(p\cdot b_{02})(\lambda_{03}-\lambda_{01})\\
       & = & (p  \cdot b_{12})(\lambda_{02}-\lambda_{03})(\lambda_{03}-\lambda_{01})\\
       & = & (\lambda_{01}-\lambda_{02})(\lambda_{02}-\lambda_{03})(\lambda_{03}-\lambda_{01}),
\end{eqnarray*}
which completes the proof.
\end{meinproof}
\par
A hyperboloid is said to be \emph{equilateral} if its asymptotic cone
is equilateral. Recall that when $Q(x)=\rho$ is an equation of a
hyperboloid, then $Q(x)=0$ is an equation of its asymptotic cone.

\begin{thm}\label{thm:hauptsatz}
If the four altitudes of a tetrahedron $\tetra$ are mutually skew,
then they are four generators of an equilateral hyperboloid $\cH$.
The Monge point $M$ is the center of $\cH$.
\end{thm}
\begin{meinproof}
Since the four altitudes are mutually skew, the scalars
$\lambda_{01}$, $\lambda_{02}$, and $\lambda_{03}$ are distinct (see
(\ref{eq:schnittbed_lambda})). Accordingly, the quadratic form $Q^*$
assumes a nonzero value at each vertex of the tetrahedron, whence
$Q^*$ cannot be the zero form. We already know from
(\ref{eq:spurnull}) that $\tr Q^*=0$ and from section
\ref{para:spurfrei} that $\rank Q^*\geq 2$. We claim that $\rank
Q^*=3$.
 \par
Assume to the contrary that $\rank Q^*=2$. Then $Q^*(x)=0$ would
define a pair of orthogonal planes and (\ref{eq:gleichung}) would be
the equation of a hyperbolic cylinder. However, the lines on a
cylinder are parallel to each other, an absurdity.
 \par
Thus $\rank Q^*=3$ and, since $Q^*(v)=Q^*(-v)$ for all $v$ in $\EE$,
equation (\ref{eq:gleichung}) describes a quadric for which the Monge
point is a center of symmetry, i.e., the reflection at the Monge
point leaves the quadric surface invariant. Therefore it is an
equilateral hyperboloid $\cH$ whose asymptotic cone has equation
$Q^*(x)=0$ and whose center is $M$. (It cannot be a hyperbolic
paraboloid, because a paraboloid does not have a center of symmetry.)
\end{meinproof}
\par
Under the assumptions of Theorem \ref{thm:hauptsatz}, we have $\tr
Q_{ijkl}\neq 0$ for every choice of indices. In this case $\cT=\RR
Q^*$ is a one-dimensional subspace of $\cQ$. Also, since there are
lines on the quadric surface $\cH$, it is a hyperboloid of one sheet.
The lines on such a quadric fall into two classes called
\emph{reguli}. Any two lines of the same regulus are skew, whereas
each line $g$ of either regulus meets all lines but one of the other
regulus; this exceptional line is parallel to $g$. Hence the four
altitudes are four lines belonging to the same regulus. By Theorem
\ref{thm:1}, each orthocentric perpendicular has three points in
common with $\cH$, whence it is a line on this surface. This ensures
that the four orthocentric perpendiculars belong to the other
regulus, i.e., not to the one containing the altitudes. This is
illustrated in Figure \ref{abb10}; however, only one orthocentric
perpendicular is actually drawn. (On a hyperbolic paraboloid two
distinct parallel lines do not exist. This demonstrates again that
the quadric surface of the altitudes cannot be a hyperbolic
paraboloid.)
 %%%%%%%%%%%%%%%%%%%%%%%%%%%%%%%%%%%%%%%%%%%%%%%%%%%%%%%%%%%%%%%%%%%%%%%%%%%%%%%
{\unitlength1.3cm%% FAKTOR!
%%   %\begin{figure}[ht]
      \begin{center}
      %%Figur Hoehen auf Quadrik
      %% 1200dpi =  472,4409448819pixel/cm
      %% 3600x3600 = 3,81 x 3,81 cm bei 2400dpi
      \begin{minipage}[t]{6.253cm}%%{4.81cm}
         \begin{picture}(4.81,3.81)
         \put(0.5,0.0)
         {\includegraphics[width=4.953cm]{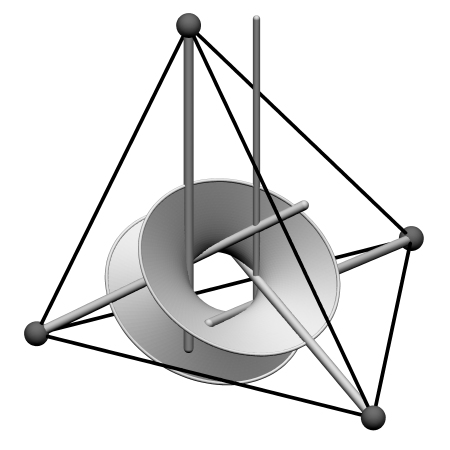}}
         \put(0.2 ,0.9){$A_0$}
         \put(3.8 ,0.1){$A_1$}
         \put(4.2 ,1.7){$A_2$}
         \put(1.5 ,3.4){$A_3$}
         \put(2.75,3.4){$n_3$}
         \end{picture}
         {\refstepcounter{abbildung}\label{abb10}
          \centerline{Figure \ref{abb10}.}}
      \end{minipage}
      \hspace{0cm}
      %%Figur Zerfallende Quadrik
      %% 1200dpi =  472,4409448819pixel/cm
      %% 3400x3000 = 3,60 x 3,18 cm bei 2400dpi
      \begin{minipage}[t]{5.98cm}%%{4.60cm}
         \begin{picture}(4.60,3.18)
         \put(0.5,0.0)
         {\includegraphics[width=4.68cm]{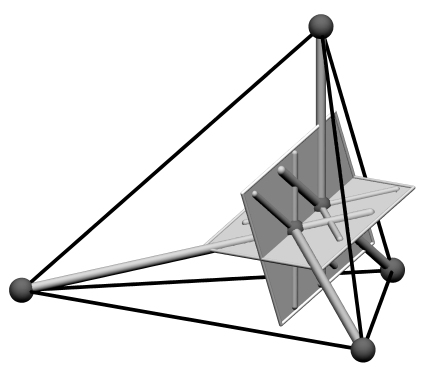}}
         \put(0.1 ,0.6){$A_0$}
         \put(3.7 ,0.1){$A_1$}
         \put(4.0 ,0.85){$A_2$}
         \put(3.4 ,2.8){$A_3$}
         \end{picture}
         {\refstepcounter{abbildung}\label{abb11}
          \centerline{Figure \ref{abb11}.}}
         \end{minipage}
      \end{center}
   %\end{figure}
}%
%%%%%%%%%%%%%%%%%%%%%%%%%%%%%%%%%%%%%%%%%%%%%%%%%%%%%%%%%%%%%%%%%%%%%%%%%%%%%%%
Next, we consider the case that two altitudes meet. Let, for example,
$\lambda_{02}=\lambda_{03}$. Then (\ref{eq:Q^*_asymm}) shows that
(\ref{eq:Q^*_def}) simplifies to
\begin{equation}\label{eq:Q^*_singulaer}
  Q^*=(\lambda_{01}-\lambda_{03})Q_{0123}.
\end{equation}
If $\lambda_{01}\neq \lambda_{03}$, then $\tr Q_{0231}\neq 0$ and
likewise $\tr Q_{0312}\neq 0$. We infer that $\cT=\RR Q^*=\RR
Q_{0123}$ is one-dimensional. Also, there is no orthocenter, but the
tetrahedron $\tetra$ is semi-orthocentric. Furthermore,
(\ref{eq:gleichung}) is the equation of a degenerate quadric surface
formed by the orthogonal midplanes $x\cdot b_{01}=0$ and $x\cdot
b_{23}=0$. The first plane carries the altitudes $h_2$ and $h_3$ and
the orthocentric perpendiculars $n_2$ and $n_3$, whereas the second
plane carries the lines $h_0$, $h_1$, $n_0$, and $n_1$. The Monge
point is the midpoint of the intersection points $h_0\cap h_1=n_2\cap
n_3$ and $h_2\cap h_3=n_0\cap n_1$. The mutual position of these
lines is portrayed in Figure \ref{abb11}.
 \par
Finally, suppose that all altitudes are concurrent. Therefore
$\lambda_{01}=\lambda_{02}=\lambda_{03}$ and all quadratic forms
$Q_{ijkl}$ are traceless. Then $\cT=\cS$ is two-dimensional. In light
of (\ref{eq:Q^*_singulaer}), $Q^*$ is a trivial quadratic form that
does not deserve any attention. It is easily seen that the solution
set of $Q_{ijkl}(x)=0$ contains all four altitudes, whence every
equation $Q(x)=0$, where $Q$ in $\cT=\cS$ is nonzero, includes all
four altitudes in its solution locus.

\para\textbf{FINAL REMARKS AND FURTHER READING.}\label{para:ende}
If we are given an equilateral hyperboloid $\cH$ of one sheet then
any plane that is perpendicular to some generator $g$ of $\cH$ meets
the surface along an equilateral hyperbola or, if it is a tangent
plane, along two orthogonal generators. This fact follows easily from
Theorem \ref{thm:gls-kegel} and is depicted in Figure \ref{abb12}.
\par
%%%%%%%%%%%%%%%%%%%%%%%%%%%%%%%%%%%%%%%%%%%%%%%%%%%%%%%%%%%%%%%%%%%%%%%%%%%%%%%
{\unitlength1cm
%%   %\begin{figure}[ht]
      \begin{center}
      %%Figur Schnitte orthogonal zur Erzeugenden
      %% 1200dpi =  472,4409448819pixel/cm
      %% 4000 x 3400 = 4,23 x 3,60 cm bei 2400dpi
      \begin{minipage}[t]{5.23cm}
         \begin{picture}(5.23,3.60)
         \put(0.5,0.0)
         {\includegraphics[width=4.23cm]{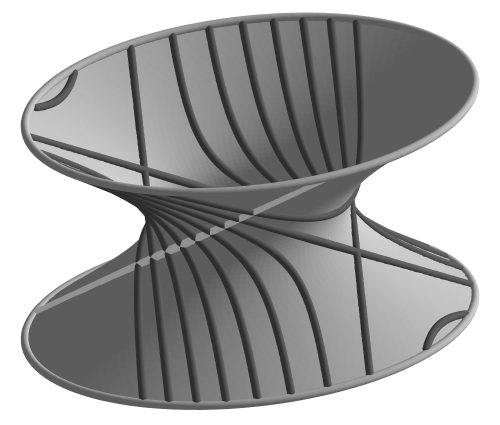}}
         \put(0.7 ,1.2){$g$}
         \end{picture}
         {\refstepcounter{abbildung}\label{abb12}
          \centerline{Figure \ref{abb12}.}}
         \end{minipage}
      \end{center}
   %\end{figure}
}%
%%%%%%%%%%%%%%%%%%%%%%%%%%%%%%%%%%%%%%%%%%%%%%%%%%%%%%%%%%%%%%%%%%%%%%%%%%%%%%%
We remark that the hyperboloid of the altitudes (which is described
in Theorem \ref{thm:hauptsatz}) meets each plane of the tetrahedron
$\tetra$ along an equilateral hyperbola. This is the result that a
conic passing through a triangle and its orthocenter is an
equilateral hyperbola \cite[p.\ 172]{semp+k-98}. The quadric surface
of the altitudes has a lot of further interesting properties (see
\cite{court-48}).
 \par
Other recent publications dealing with the altitudes of a
tetrahedron, the altitudes of simplexes in higher dimensional spaces,
and other topics that are related to this circle of ideas are:
 \cite[p.~371]{bolt+jagl-69},
 \cite{coud+b-35},
 \cite{court-34}, \cite{court-37}, \cite{court-42},
 \cite{court-46}, \cite{court-48}, \cite{court-53},
 \cite{court-64},
 \cite{fritsch-76},
 \cite{gerb-73},
 \cite{mandan-58a}, \cite{mandan-58b}, \cite{mandan-58c},
 \cite{mandan-61}, \cite{mandan-62a}, \cite{mandan-62b},
 \cite{mandan-65},
 \cite{micu-88},
 \cite[p.~376]{rose+jagl-69},
 \cite{satya-87}.

%%\bibliographystyle{plain}
%%\bibliography{D:/forschung/separata/hoehen}

%%ACHTUNG BIBTEX sortiert die Mandan-Arbeiten schlecht!!!

\noindent
Hans Havlicek\\
\emph{Institut f\"ur Geometrie, Technische Universit\"at, Wiedner Hauptstra{\ss}e
8--10, A--1040 Wien, Austria}
\\ \emph{havlicek@geometrie.tuwien.ac.at}

 \vspace{1ex}
 \noindent
Gunter Weiss\\
\emph{Institut f\"ur Geometrie, Technische Universit\"at, Zellescher Weg
12--14, D--01062 Dresden, Germany
\\ weiss@math.tu-dresden.de}

\end{document}